# Fatigue reliability analysis of offshore wind turbines under combined wind-wave excitation via DPIM

Jingyi Ding[a], Hanshu Chen[a,*], Xiaoting Liu[b], Youssef F. Rashed[c], Zhuojia Fu[a,*]

[a] Center for Numerical Simulation Software in Engineering and Sciences, College of Mechanics and Engineering Science, Hohai University, Nanjing, Jiangsu, 211100, China
[b] Institute of Science and Technology Research, China Three Gorges Corporation, Beijing 210098, China
[c] Department of Structural Engineering, Cairo University, Giza, Egypt

**Abstract:** As offshore wind turbines develop into deepwater operations, accurately quantifying the impact of stochastic excitations in complex sea environments on offshore wind turbines and conducting structural fatigue reliability analysis has become challenging. In this paper, based on long-term wind-wave reanalysis data from a site in the South China Sea, a novel direct probability integral method (DPIM) is developed for the stochastic response and fatigue reliability analyses of the key components for the floating offshore wind turbine structures under combined wind-wave excitation. A 5MW floating offshore wind turbine is considered as the research object, and a fully coupled dynamic response analysis of the wind turbine system is conducted to calculate the short-term fatigue damage value of tower base and blade root. The DPIM is applied to calculate the fatigue reliability of the wind turbine structure. The accuracy and efficiency of the proposed method are validated by comparing the obtained results with those of Monte Carlo simulations. Furthermore, the results indicate that the fatigue life of floating offshore wind turbine structures under combined wind-wave excitation meets the design requirements. Notably, the fatigue reliability of the wind turbine under aligned wind-wave condition is lower compared to misaligned wind-wave condition.

## 1 Introduction

Offshore wind energy, with its vast resource reserves and flexible environmental constraints, has become an important development direction for future wind energy [1]. However, as offshore wind power projects gradually advance into deepwater operations, the wind turbines will be subjected to complex, variable, and alternating random loads, which can easily induce cumulative fatigue damage and subsequently threaten the overall safety and stability of the structure [2, 3]. Therefore, conducting a fatigue reliability analysis of the offshore wind turbines under stochastic environmental excitations is of significant importance.

In recent years, researchers have conducted extensive studies on the reliability analysis of offshore wind turbine units. For instance, Colone et al. [4] examined the impact of turbulence and wave loads on the fatigue reliability of offshore wind turbine pile foundations, emphasizing the need for accurate environmental modeling to improve cost-effective and reliable monopile design. Wilkie et al. [5] developed a computational framework using Gaussian process regression to effectively assess the fatigue reliability of offshore wind turbine substructures, providing practical insights for design and reliability in European waters. Fu [6] investigates the fatigue reliability of the wind turbine tower flange and bolt under random wind loads, proposing a fatigue probability calculation method based on probability density evolution for accurate and quantitative prediction of structural reliability under fluctuating conditions. Zhao et al. [7] proposed a fatigue reliability analysis method using a surrogate model, C-vine copula, and Monte Carlo simulation, applied to assess fatigue reliability at three critical locations on a floating offshore wind turbine (FOWT). Although there has been considerable research on the reliability analysis of wind turbines under the combined wind-wave excitations have been studied, the research on the fatigue reliability analysis of wind turbine structures considering multiple uncertainty factors remains limited.

Currently, structural reliability analysis methods mainly include the first-order second-moment method (FORM), second-order second moment methods (SORM), random sampling simulation methods [8], and probability density evolution methods [9, 10]. FORM and SORM are commonly used for static analysis and are not suitable for the dynamic systems of wind turbines. Random sampling simulation methods, such as Monte Carlo Simulation (MCS), subset simulation, and importance sampling, yield accurate results but come with high computational costs. Li and Chen [9] established the generalized density evolution equation (GDEE) based on the principle of probability conservation and developed the probability density evolution method. However, analytical solutions for GDEE are challenging for most engineering problems. To address this issue, this paper extends a novel direct probability integral method (DPIM) [11, 12] for the stochastic response analysis of wind turbine structures subjected to stochastic environmental excitations, aiming to achieve efficient fatigue reliability analysis of FOWT under combined wind-wave excitation.

This study focuses on the 5 MW spar wind turbine developed by the National Renewable Energy Laboratory (NREL) [13]. An efficient assessment method for the fatigue reliability analysis of the key components for wind turbine structures with multiple uncertainty factors is further proposed based on DPIM. This model incorporates the long-term joint probability distribution of wind and wave conditions relevant to the specific maritime area for time-domain analysis. Using the rainflow counting method, the stress spectra of danger nodes in the support structure are extracted. Then, the number of cycles is computed via S-N curves [14], and Palmgren-Miner's cumulative fatigue damage theory [15] is applied to estimate short-term fatigue damage values. Finally, the fatigue reliability analysis of the wind turbine structure using DPIM is achieved by the proposed method.

The paper is organized as follows: Section 2 presents a detailed overview of the FOWT model and the fatigue analysis method via DPIM used in this study. Section 3 demonstrates the effectiveness of the proposed method. Moreover, the influence of different wind-wave conditions on stochastic response, fatigue damage, and fatigue reliability is further investigated. Lastly, Section 4 emphasizes the main conclusions and insights drawn from this research.

## 2 FOWT model and fatigue reliability analysis method

### 2.1 Numerical model of FOWT

The fatigue analysis is conducted for an NREL 5MW OWT, the platform of which is the OC3-Hywind spar-buoy [16] (Fig. 1). The main properties of the wind turbine are listed in Table 1. The height of the wind turbine, the top of the tower, and the base of the tower are 90 m, 87.6 m, and 10 m above the SWL (still water level), respectively. The draft of the platform is 120 m and the CM (center of mass) of the floating platform (including the ballast) is 89.92 m below SWL. The floating platform is moored by three equally spaced catenary lines spread symmetrically about the Z-axis, one of which is directed along the negative direction of the X-axis. The fairleads are located at a depth of 70 m below the SWL and the depth of anchors is 320 m below the SWL.

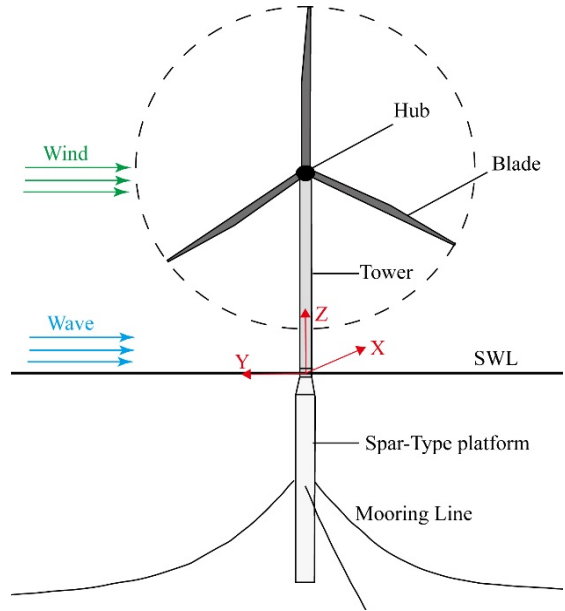

Fig. 1 5MW OC3 Spar-type FOWT

Table 1 The main properties of the 5MW OC3 Spar-type FOWT

| Parameter | Value |
| --- | --- |
| Rating | 5 MW |
| Rotor orientation, configuration | Upwind, 3 blades |
| Rotor and hub diameter, hub height | 126 m, 3 m, 90 m |
| Cut-in, rated, cut-out wind speed | 3 m/s, 11.4m/s, 25m/s |
| Cut-in, rated root speed | 6.9 rpm, 12.1 rpm |
| Elevation to tower base above SWL | 10 m |
| Tower base diameter, tower base thickness | 6.5 m, 0.027 m |

In this study, a coupled areo-hydro-servo-elastic tool OpenFAST [17], which is developed and verified by NREL, is used to obtain the dynamic responses of critical components of FOWT under different environmental conditions. The aerodynamic loads on the blades are calculated based on BEM (Blade Element Momentum) theory in the module of AeroDyn, which processes turbulent wind data produced by TurbSim. The Kaimal spectrum is used to simulate the turbulent incoming wind field for FOWT. The wind loads on the tower are calculated based on potential flow theory. For the hydrodynamic loads, both potential theory and Morison's equation are applied through the HydroDyn module. The viscous drag forces are accounted for by considering the drag term in Morison's equation. The second-order wave forces are neglected, as they are minimal compared to the first-order force acting on the spar-type floater [18]. The motion equation for FOWTs is derived using Kane's dynamics.

## 2.2 Probabilistic modeling of long-term joint wind and wave loads

FOWTs operate in variable natural environments for extended periods, facing multiple random environmental factors primarily driven by wind and wave loads. Following IEC 61400-3 [19], this study selects three environmental parameters that influence wind and wave loads, mean wind speed ($V_W$), significant wave height ($H_S$), and spectral peak period ($T_P$). The wind speed can define the characteristics of the wind states, the others can determine the characteristics of the sea states. Other environmental conditions, such as wind and wave directions, can also significantly impact fatigue reliability analysis. In this study, the fixed wind direction is aligned with the positive X-axis. The impact of varying wind and wave directions is achieved by adjusting the angle of wave inflow. The wave angles are set to 0°, 30°, 60°, and 90°, respectively.

Applying reanalysis data from the South China Sea station (21°N, 113°E) from 1979 to 2016 [20], and then simulating the environmental data for this sea area based on the joint distribution model of combined wind-wave excitation established by Song et al [21]. The probability distribution models of the environmental random variables are presented in Table 2. The corresponding PDF curves of environmental random variables are displayed in Fig. 2.

Table 2 Probability distribution models of the environmental random variables

| Variable | Distribution type | PDF $f(x)$ | Distribution parameters |
|---|---|---|---|
| $V_W$/(m·s$^{-1}$) | Truncated Weibull | $\frac{a}{b}\left(\frac{x}{a}\right)^{b-1}\exp\left(-\left(\frac{x}{a}\right)^b\right), x\in[3,25]$ | $a=11.9799$<br>$b=2.8005$ |
| $H_S$/m | Lognormal | $\frac{1}{x\sigma\sqrt{2\pi}}\exp\left(\frac{-(\ln x-\mu)^2}{2\sigma^2}\right)$ | $\mu=0.4887$<br>$\sigma=0.4489$ |
| $T_P$/s | Lognormal | $\frac{1}{x\sigma\sqrt{2\pi}}\exp\left(\frac{-(\ln x-\mu)^2}{2\sigma^2}\right)$ | $\mu=2.0759$<br>$\sigma=0.1547$ |

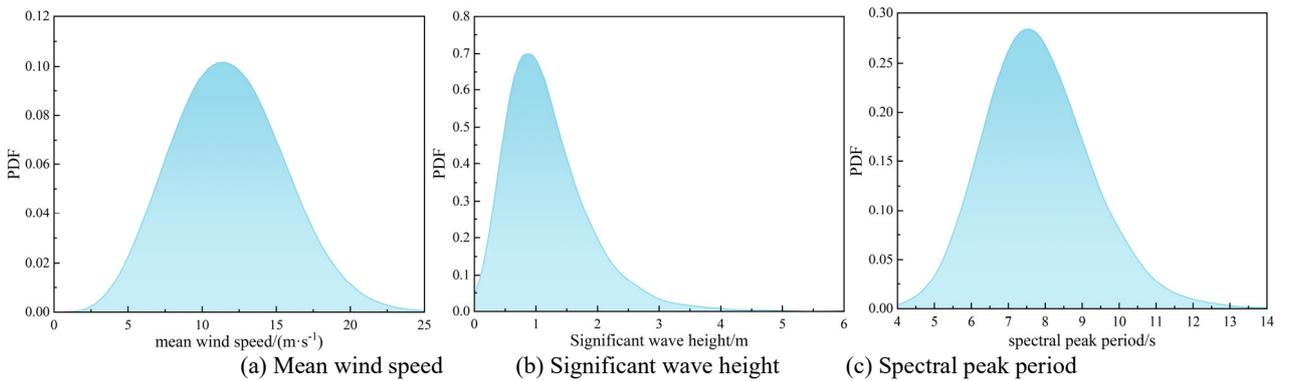

(a) Mean wind speed    (b) Significant wave height    (c) Spectral peak period
Fig. 2 PDF curves of wind and wave loading parameters under DPIM

## 2.3 Fatigue damage location and stress calculation

In this study, two hot-spot locations of the FOWT system, including the tower base and the blade root, are selected to evaluate fatigue reliability. A significant aspect of fatigue damage in the time domain is identifying the stress time series. Therefore, obtaining the stress time series from dynamic simulations at these locations is essential.

The base section of the tower is simplified as a thin-walled cylindrical structure subjected to axial and shear stresses. To simplify the calculation, the effects of bolts and the connection components between the tower and platform are ignored, and the cross-section of the tower is treated as plane stress that varies along the circumference. The fatigue damage caused by shear stress is considerably lower than that caused by axial stress [22], so only axial stresses are used to calculate fatigue damage in this study. The axial stress at the tower base is calculated at 7 points around its circumference (Fig. 3). Based on the small deformation assumption, the nominal axial stress can be calculated using Eq.(1).

$$\sigma = \frac{N_z}{A} + \frac{M_y}{I_y} \cdot r \cdot \cos\alpha - \frac{M_x}{I_x} \cdot r \cdot \sin\alpha \qquad (1)$$

where $N_z$ is axial force, $A$ is the nominal cross sectional area, $M_x$ and $M_y$ are the tower base roll and pitching moment, respectively, $I_x$ and $I_y$ are the sectional moments of the area to X-axis and the Y-axis, respectively, and $\alpha$ is the angle measured in the counterclockwise direction from the negative X-axis to the calculated point.

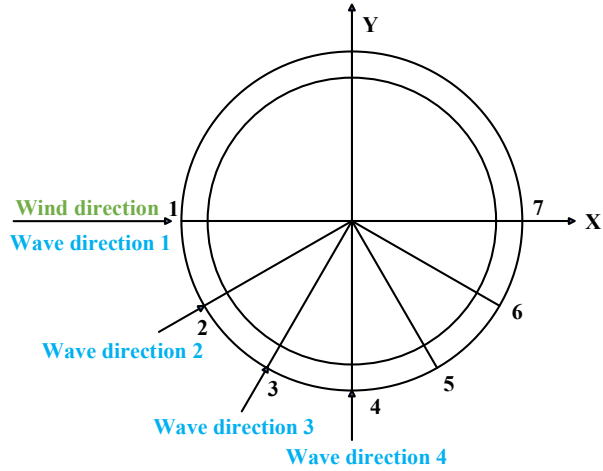

Fig. 3 Top view of the tower base

In contrast, the stress states at the blade root are more complex. The stress at the blade root is also composed of axial stress and shear stress, both of which are considered in this study. The wind turbine blade undergoes bending moments and shear forces in three directions as it rotates continuously. A fixed coordinate system is established to accurately describe these loads on the blade (Fig. 4). This coordinate system remains stationary and does not change with the rotation of the blade, allowing for consistent analysis of the forces acting on the blade. To simplify the calculation, the blade root cross-section is treated as a circular ring-shaped area. The axial stress and shear stress can be calculated using Eq.(2) and Eq.(3) respectively, which is

$$\sigma_0 = \frac{\sqrt{M_x^2 + M_y^2}}{W_n} + \frac{F_z}{A} \qquad (2)$$

$$\tau_0 = \frac{\sqrt{F_x^2 + F_y^2}}{A} + \frac{M_z}{W_p} \qquad (3)$$

where $M_x$, $M_y$ and $M_z$ are the in-plane, out-of-plane, and pitching moment, respectively; $F_x$, $F_y$ and $F_z$ are the out-of-plane shear, in-plane shear, and axial force, respectively; $A$ is the nominal cross-sectional area; $W_n$ and $W_p$ are the section modulus and section modulus in torsion, respectively.

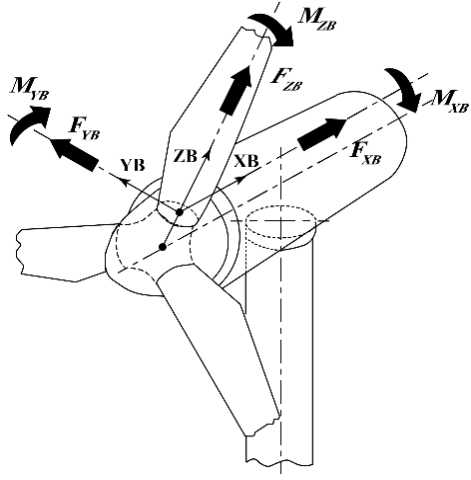

Fig. 4 The coordinate system of wind turbine blade

## 2.4 Rainflow counting method

The rainflow counting method was introduced by Matsuishi and Endo in 1968 [23] and has since become a widely accepted method for analyzing random signals in the context of fatigue analysis. The primary factors influencing the fatigue life of structures include stress amplitude, the number of cycles, and mean stress amplitude. The rain flow counting method, a statistical approach commonly employed to predict fatigue life and evaluate fatigue damage, is often used to quantify these parameters.

In this paper, the time domain simulation of the 5MW spar offshore wind turbine is first performed. Then the rain-flow counting method is employed to count the stress cycles in the stress time history. To account for the impact of the mean stress value on the structure's fatigue life, the Goodman method (Eq.(4)) is applied to correct the average stress.

$$\sigma_i^{RF} = \sigma_i^{R} \cdot \left( \frac{\sigma^{ult} - |\sigma^{MF}|}{\sigma^{ult} - |\sigma_i^{M}|} \right)^{\varepsilon} \tag{4}$$

where $\sigma_i^{RF}$ is the cycle's stress range about a fixed load-mean value; $\sigma_i^{R}$ is the i$^{th}$ cycle's range about a load mean stress of $\sigma_i^{M}$; $\varepsilon$ is Goodman exponent. In this work, it is defined as 1.0. $\sigma^{ult}$ is the highest stress of the cross-section (in absolute value) before failure based on the ultimate strength. $\sigma^{MF}$ is the fixed mean stress value based on the time series of stress. To eliminate the mean stress effect, $\sigma^{MF}$ can be set to zero.

## 2.5 S-N curve

The S-N curve defines a material's fatigue resistance by illustrating the maximum number of cycles it can withstand at a specific stress amplitude without experiencing fatigue damage. According to the DNV standard, the basic formula for calculating the S-N curve under different environmental conditions is as follows:

$$\lg N = \lg \bar{a} - m \lg \Delta\sigma \tag{5}$$

Where $N$ is the predicted number of cycles to failure for the stress range $\Delta\sigma$; $\bar{a}$ and $m$ are the slope and the intercept parameters of the S-N curve. The blades are made of composite materials, while the tower is made of steel, $m$ is defined as 8 and 3, respectively.

## 2.6 Fatigue cumulative damage theory

Fatigue damage refers to the extent of material degradation due to fatigue loading, often represented by the dimensionless parameter $D$. When $D = 0$, it indicates that the material has not undergone any fatigue damage. Conversely, when $D > 1$ suggests that the material has exhausted its fatigue life. In this study, P-M linear fatigue cumulative damage theory is used to evaluate fatigue damage. The foundational formula is as follows:

$$D_j^{ST} = \sum_i \frac{n_{ji}}{N_{ji}} \tag{6}$$

$$DR_j^{ST} = \frac{D_j^{ST}}{T_j} \tag{7}$$

$$D = \sum_{j=1}^{N} DR_j^{ST} P_q T \tag{8}$$

Where $D_j^{ST}$ is the fatigue damage value caused by ith stress cycle; $n_{ji}$ is the number of cycles for the i$^{th}$ stress cycle obtained by the rainflow counting method; $N_{ji}$ is the number of failures under the i$^{th}$ stress cycle; $T_j$ is the simulation time ($T_j = 600$ s in this paper); $N$ indicates the number of load cases for fatigue damage calculation; $T$ stands for the total operating time of the wind turbine and $D$ represents the cumulative fatigue damage.

## 2.7 Fatigue reliability analysis via DPIM

Assume the stochastic sources for the multi-degree-of-freedom (MDOF) nonlinear system of floating wind turbines under random external excitation (such as wind, waves, etc.) all come from $\boldsymbol{\Theta}$, the principle of probability conservation [9] for the key components of wind turbines is expressed as:

$$\int_{\Omega_Y} p_Y(\mathbf{y},t) \mathrm{d}\mathbf{y} = \int_{\Omega_\Theta} p_\Theta(\boldsymbol{\theta},t) \mathrm{d}\boldsymbol{\theta} \tag{9}$$

in which $p_\Theta(\boldsymbol{\theta}, t)$ and $p_Y(\mathbf{y}, t)$ indicate PDF of input random variables $\boldsymbol{\theta}$ and output vector $\mathbf{y}$, respectively. The corresponding relationship between the input random variables $\boldsymbol{\theta}$ and output vector $\mathbf{y}$ can be described by a deterministic mapping $\mathrm{G}$:

$$\mathrm{G}: \mathbf{Y}(t) = \mathbf{g}(\Theta, t) \tag{10}$$

To further explicitly characterize the uncertainty propagation from input variables $\boldsymbol{\theta}$ into output response vector $\mathbf{y}$, based on the deterministic mapping G and Dirac delta function, PDIE of MDOF nonlinear system of floating wind turbines is derived, i.e.,

$$p_\mathbf{Y}(\mathbf{y},t)=\int_{-\infty}^{\infty}\cdots\int_{-\infty}^{\infty}p_\mathbf{\Theta}(\mathbf{\theta})\delta[\mathbf{y}-\mathbf{g}(\mathbf{\theta},t)]\,d\mathbf{\theta} \qquad (11)$$

Since in practice it is not necessary to obtain a joint PDF for all responses, only a few or even a single response is of primary concern. The dimension reduction for MDOF nonlinear system of floating wind turbines is accomplished using the property of Dirac delta function and implementing the marginal integral of Eq. (11). As a result, PDF of concerned response $y_\ell(t)$ is derived, which is

$$p_{Y_1}(y_1,t)=\int_{-\infty}^{\infty}L\int_{-\infty}^{\infty}\delta[y_1-g_1(\mathbf{\theta},t)]p_\mathbf{\Theta}(\mathbf{\theta})\,d\mathbf{\theta} \qquad (12)$$

To address the issues caused by the singularity of the Dirac function in solving Eq. (12), Chen et al. [10] introduced the techniques of the partition of probability space and smoothing of the Dirac delta. Consequently, the numerical computation formula for Eq. (12) is obtained as follows:

$$p_{Y_\ell}(y_\ell,t)=\sum_{q=1}^{N}\left\{\frac{1}{\sqrt{2\pi}\sigma}e^{-[y_\ell-g_\ell(\mathbf{\theta}_q,t)]^2/2\sigma^2}P_q\right\} \qquad (13)$$

where $N$ indicates the number of representative points generated by the techniques of the partition of probability space, $g_\ell(\mathbf{\theta}_q,t)$ means the representative response of $q$-th representative points, $P_q$ is the assigned probability of $q$-th representative points, and $\sigma$ denotes the smoothing parameter.

Then, for achieving the fatigue reliability analysis of key components for the floating wind turbines via DPIM, the dynamic performance function is denoted as

$$M:Z(t)=F(\mathbf{\Theta},t) \qquad (14)$$

in which $Y(\mathbf{\Theta},t)$ is expressed as

$$Z(t)=B-g_{\text{ext}}(\mathbf{\Theta},t) \qquad (15)$$

where $B$ denotes the threshold, $g_{\text{ext}}(\mathbf{\Theta},t)$ means the equivalent extreme value of fatigue response of wind turbines. Furthermore, the fatigue reliability function of floating wind turbines is formulated by

$$\begin{aligned}R_s(t)&=\Pr\{Z(\tau)\in\Omega_{Z,s},\tau\in(0,t]\}\\&=\Pr[Z>0]=\int_0^{\infty}p_Z(z,t)\,dz\end{aligned} \qquad (16)$$

where $\Omega_{Z,s}$ is safe domain of stochastic response $Z(t)$, and $p_Z(z,t)$ denotes the PDF of $Z(t)$. According to PDIE, the fatigue reliability function in Eq. (16) can be evaluated by

$$\begin{aligned}R_s(t)&=\Pr[Z>0]=\int_0^{\infty}p_Z(z,t)\,dz\\&=\int_0^{\infty}\int_{\Omega_\mathbf{\Theta}}\delta[z-F(\mathbf{\theta},t)]p_\mathbf{\Theta}(\mathbf{\theta})\,d\mathbf{\theta}\,dz\end{aligned} \qquad (17)$$

Introducing Heaviside function H, the formulation for calculating the fatigue reliability of floating wind turbines in the framework of DPIM is given by

$$\begin{aligned}R(t)&=\int_{-\infty}^{\infty}H[F(\mathbf{\theta},t)]p_\mathbf{\Theta}(\mathbf{\theta})\,d\mathbf{\theta}\\&=\sum_{q=1}^{N}\{H[F(\mathbf{\theta}_q,t)]P_q\}\end{aligned} \qquad (18)$$

Equation (18) introduces Heaviside function to avoid the smoothing process of Dirac delta function, and the numerical solution of Eq. (18) only requires the application of the techniques of the partition of probability space, which is easier to solve.

# 3  Numerical examples

In this paper, a two-step procedure for constructing representative points based on GF-discrepancy is employed, and the number of the representative points is set as 1000. As a result, the stochastic and fatigue reliability analyses of tower base and blade root for the 5 MW spar wind turbine are achieved based on DPIM. Furthermore, MCS method results are used as a benchmark solution, the sample size for the method is set to 10,000.

## 3.1  Stochastic response analysis of wind turbines under combined wind-wave excitation

The comparison of the average stress at key points at the tower base under different wind and wave angle conditions (0°, 30°, 60°, and 90°) is illustrated in Fig. 5. By comparing the mean axial stress values at various points of the tower base, it is evident that the maximum mean axial stress occurs at Node 7. This suggests a potential fatigue risk that warrants further investigation. Additionally, according to Eq. (1), the axial stress varies with an angle like a sinusoidal function. This indicates that as a point moves further away from the direction of the wind and waves in the Y direction, both the absolute value of the mean axial stress and the range of axial stress decrease. Furthermore, the absolute mean axial stress values on the lee side are slightly higher than those facing the waves. This increase is due to the superposition of the axial stresses resulting from the flatwise bending moment and the axial force acting on the points located on the lee side. Additionally, despite the changes in the wind and wave angles, the difference in average stresses at various nodes is not significant. Thus, the node stresses are not sensitive to variations in wind and wave angles.

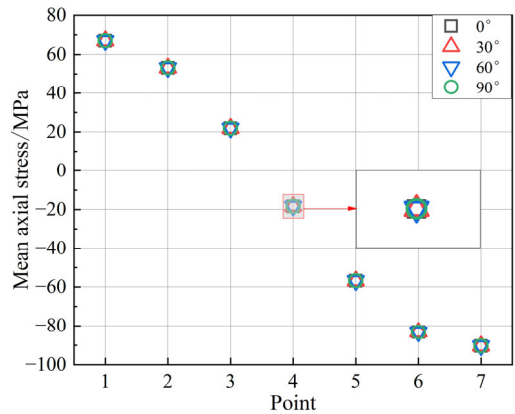

Fig. 5 Comparison of mean axial stress under different wind and wave angle conditions

The average stress values at different points at the tower base under the aligned wind-wave condition, calculated using the DPIM and MCS methods, are shown in Fig. 6 (a), respectively. Specifically, the PDF curves of the stress at point 7 of the tower base and the blade root are shown in Fig. 6 (b) and Fig. 6 (c), respectively. The comparative results show that the absolute stress at Node 7 ranges primarily from 80 MPa to 120 MPa. In contrast, the stress at the blade root typically falls between 20 MPa and 45 MPa. This indicates that the average stress on the supporting structure is higher than that at the blade root. Additionally, the results from both DPIM and MCS are consistent with each other.

Furthermore, regarding computational efficiency, the CPU time for DPIM is approximately 22,910

s, while the MCS method takes 611,770 s. This means that DPIM is over 20 times more efficient, demonstrating the accuracy and high efficiency of the proposed method for uncertainty qualification analysis of large-scale wind turbine structures.

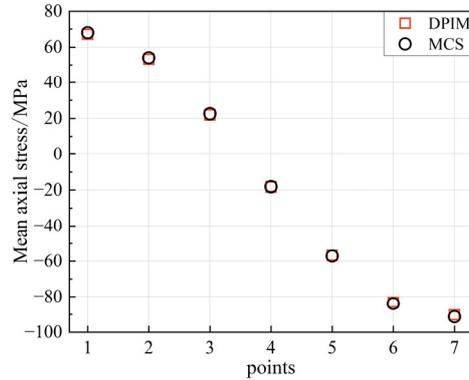
(a) mean axial stress at point 7 under different methods

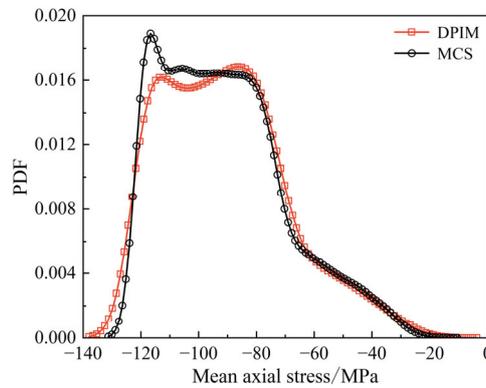
(b) PDF curve of mean axial stress at point 7 under different methods

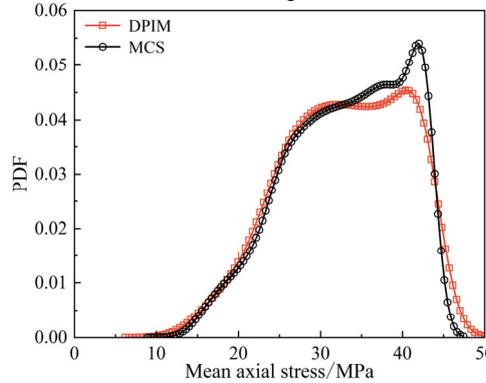
(c) PDF curve of mean axial stress at blade root under different methods
Fig. 6 Comparison of mean stress results under different methods

## 3.2 Fatigue damage analysis of wind turbines under combined wind-wave excitation

To comprehensively assess the impact of fatigue accumulation effects under long-term operating conditions, the calculated fatigue damage values are first converted to a damage rate per second, which is then multiplied by the corresponding calculation duration to quantify the accumulation of damage over time. The PDF curves of the fatigue damage at the danger point (Node 7) of tower base and blade root after 20 years of operation under different wind and wave directions are illustrated in Fig. 7. When the wind and wave angle is 90°, the probability that the fatigue damage value at the tower base and blade root is less than 1 is the highest compared to all other angle conditions. Additionally, the probability density function (PDF) curve shows peaks in the ranges of [0.5, 0.8] and [0, 0.25], indicating the most likely variation ranges of fatigue damage for the tower base and blade root, respectively.

In contrast, when the wind and waves are aligned, the range of the structural fatigue damage values increases, and the probability of values greater than 1 also rises, making the probability of being less than 1 the lowest among all angle conditions. It can be concluded that the wind and wave angles have a significant impact on the fatigue damage of wind turbine structure, with relatively lower fatigue damage observed at 90°, while the probability of fatigue damage values exceeding 1 is highest when the wind and waves are aligned. The diagram comparison indicates that the likelihood of the fatigue damage value at the blade root being less than 0.5 is significantly higher than at the tower base. This suggests that the tower base is more susceptible to fatigue failure compared to the blade root. These findings align with the average stress results discussed in Section 3.1.

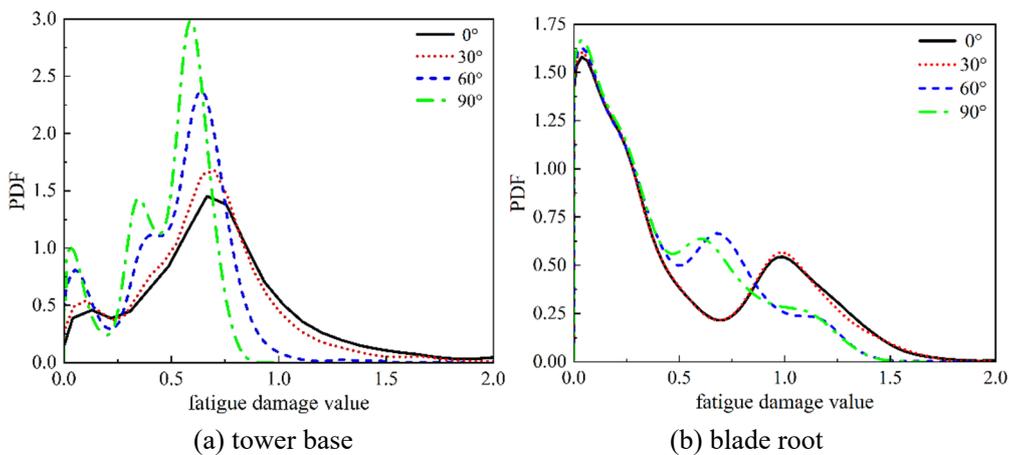

(a) tower base　　　　　　　　　　　(b) blade root

Fig. 7 Comparison of PDF curves of fatigue damage value at key components under different wind-wave angles

## 3.3 Fatigue reliability analysis of wind turbines under combined wind-wave excitation

In this section, the impact of wind and wave angle on the fatigue reliability of critical points under different operating years has been further studied. As shown in Fig. 6 (a), when the wind and waves are aligned, the assessment of fatigue reliability for the tower base of wind turbines is implemented based on fatigue damage at the danger point (point 7). As the wind and wave angle increases from 0° to 90°, the fatigue reliability of wind turbines shows a significant upward trend, which corresponds with the previously analyzed fatigue damage value results in Section 3.2.

As the operational lifespan of the wind turbine increases, the fatigue reliability of the tower base at various points shows a gradual decline. Under the aligned wind-wave condition, when the operational lifespan extends to the design life of 20 years, the fatigue reliability of the wind turbine is 0.855. This relatively high value indicates that the wind turbine structure possesses a significant level of fatigue reliability within the design life. However, when the operational lifespan is extended to 25 years, the fatigue reliability is significantly reduced to 0.704, indicating the adverse effects of long-term operation on the fatigue performance of the wind turbine structure. This result indicates that once the operational lifespan of the wind turbine exceeds its design life, the fatigue reliability of the wind turbine structure will gradually decline, with this downward trend becoming increasingly pronounced as the operational lifespan continues to increase.

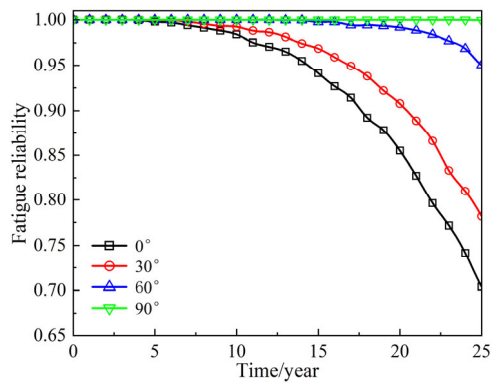

Fig. 8 Comparison of fatigue reliability of tower base under different wind-wave conditions for different operating years

Furthermore, the fatigue reliability of blade root for the wind turbines is also given in Fig. 9. It can be seen that the combined effects of wind and waves significantly impact the fatigue reliability of the blade. The fatigue reliability begins to sharply decline when the operational lifespan is much less than the designed lifespan (20 years) of wind turbines. Moreover, when the angle between the wind and waves is small (e.g., 0°and 30°,), the stochastic environmental excitations significantly affect the fatigue reliability of the wind turbine blades. In comparison to the misaligned wind-wave condition, the aligned wind-wave condition results in a greater impact of these excitations, leading to a decrease in fatigue reliability from 0.864 at a design life of 20 years to 0.746 at 25 years. Notably, as illustrated in Figs. 8 and 9, after the wind turbine has been in operation for more than 20 years, the fatigue reliability of the wind turbine structure sharply decreases. Therefore, it is essential to conduct regular inspections, maintenance, or reinforcements in high-risk areas to ensure the long-term safe and stable operation of the wind turbines.

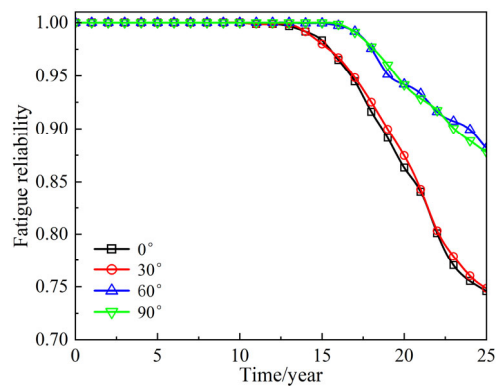

Fig. 9 Comparison of fatigue reliability of blade root under different wind-wave conditions for different operating years

## 4  Conclusions

In this paper, an efficient analysis method for the fatigue reliability of the 5 MW spar offshore wind turbine structure is proposed based on DPIM. The main conclusions are as follows:

(1) This paper develops DPIM to address the stochastic response and structural fatigue reliability analyses of key components for floating offshore wind turbines under combined wind-wave excitation. The proposed uncertainty quantification method can be applied as an alternative tool for floating offshore wind turbine structures with multiple uncertain factors.

(2) Compared to typical MCS, DPIM has high accuracy and significantly improves computational efficiency, making the assessment for the fatigue reliability analysis of key components for large floating offshore wind turbine structures more efficient.

(3) The results of the uncertainty quantification analysis of the wind turbine under the combined wind-wave excitation indicate that the fatigue reliability of the tower base and blade root significantly decreases over time. Under the aligned wind-wave condition, when the wind turbine operates up to its design life of 20 years, the fatigue reliabilities of the tower base and blade root are 0.855 and 0.864, which is lower than under other angle conditions but still demonstrates high reliability. When the operational life cycle reaches 25 years, the fatigue reliabilities of the tower base and blade root are dropped to only 0.704 and 0.746, indicating a sharp decline in reliability. At this stage, regular inspections, maintenance, or reinforcements of high-risk areas are required to ensure the long-term safe and stable operation of the wind turbine.

**Acknowledgments**

The supports of the National Natural Science Foundation of China (Grant Nos. 12372196, 12402238), Jiangsu Funding Program for Excellent Postdoctoral Talent (Grant No. 2023ZB506) and Postdoctoral Fellowship Program of CPSF (Grant No. GZC20230667) are much appreciated.